\documentclass[reqno]{amsart}
\usepackage{float}
\usepackage{amsmath}
\usepackage{graphicx}
\usepackage{latexsym}
\usepackage{amsfonts}
\usepackage{amssymb}
\usepackage{color}
\theoremstyle{plain}
\newtheorem{theorem}{Theorem}
\newtheorem{corollary}[theorem]{Corollary}
\newtheorem{lemma}[theorem]{Lemma}
\newtheorem{proposition}[theorem]{Proposition}
\theoremstyle{definition}

\newtheorem{definition}[theorem]{Definition}

\newtheorem*{remark*}{Remark}

\newcommand{\pr}{\mathbf P}
\newcommand{\e}{\mathbf E}

\begin{document}
\title[Maximum on a random time interval]{Maximum on a random time interval of a random walk with infinite mean}

\author[Denisov]{Denis Denisov}
\address{School of Mathematics, University of Manchester, Oxford Road, Manchester M13 9PL, UK}
\email{denis.denisov@manchester.ac.uk}

\begin{abstract}
    Let $\xi_1,\xi_2,\ldots$ be independent, identically distributed random variables with infinite mean $\e[|\xi_1|]=\infty.$ Consider a random walk $S_n=\xi_1+\cdots+\xi_n$, a stopping time $\tau=\min\{n\ge 1: S_n\le 0\}$ and let $M_\tau=\max_{0\le i\le \tau} S_i$. We study the asymptotics for $\pr(M_\tau>x),$ as $x\to\infty$. 
\end{abstract}
\maketitle

\section{Introduction}
Let $\xi,\xi_1$, $\xi_2$, \ldots\ be independent random variables with
a common distribution~$F$.
Consider a random walk $S_0=0$, $S_n=\xi_1+\ldots+\xi_n$ and a stopping time  
\[
    \tau:=\inf\{n\ge 1: S_n<0\}.
\]
Let $M_\tau:=\max_{0\le i\le \tau} S_i$ and $M=\sup\ \{S_n,\ n\ge 0\}$.
 We will consider the random walks with infinite or undefined mean ($\e[|\xi_1|]=\infty$) under the assumption that  $S_n\to-\infty$ a.~s. 
It is well known that the latter assumption is equivalent to 
$M<\infty$ a.~s. and to $\e[\tau]<\infty$ (see Theorem 1 in~\cite[Chapter XII, Section 2]{F}). 

In the infinite-mean case, an important role is played by {\it the negative truncated mean function}
\begin{eqnarray*}
m(x) &\equiv& {\bf E}\min\{\xi^-,x\}
=\int_0^x {\bf P}\{\xi^->y\}\,dy,\quad x\ge0,
\end{eqnarray*}
where $\xi^-=\max\{-\xi,0\}$; the function $m(x)$ is continuous, increasing,
$m(0)=0$ and $m(x)>0$ for any $x>0$.
It is known that if  $\e |\xi|=\infty$, then $S_n\to-\infty$ 
a.s. as $n\to\infty$ if and only if
\begin{equation}\label{finite.M}
    K:=\int_0^\infty\frac{x}{m(x)}\,F(dx)\ \mbox{ is finite,}
\end{equation}
see Corollary 1 in~\cite{E}. 

The aim of this paper is to study the asymptotics for $\pr(M_\tau>x)$ in the infinite-mean case. 
In the finite-mean case, under the assumption that $F\in \mathcal S^*$
it was shown  by Asmussen~\cite{Asmussen1998}, see also~\cite{HRS1997} for the regularly varying case,  that 
\begin{equation}\label{asym.Mtau}
\mathbf P(M_\tau>x)\sim \mathbf E\tau \overline F(x),
\end{equation}
where $\overline F(x) = 1-F(x) = \pr(\xi_1>x).$
The class  $\mathcal S^*$ of strongly subexponential distributions was introduced by Kl\"uppelberg~\cite{Kl} and is defined as follows, 
\begin{definition}\label{sstar}
    A distribution function $F$ with  finite $\mu_+=\int_0^\infty \overline F(y)dy<\infty$ belongs to the class $\mathcal S^*$ of strong subexponential distributions  if 
 $\overline F(x)>0$ 
for all $x$ and 
\[
\frac{\int_0^x \overline F(x-y)\overline F(y)dy}{\overline F(x)}
\to 2 \mu_+,\quad \mbox{ as } x\to \infty
\]
\end{definition}
This class is a proper subclass of class $\mathcal S$ of subexponential distributions. It is shown 
in~\cite{Kl} that the 
Pareto, lognormal and Weibull distributions belong to the class $\mathcal S^*$
as well. 

The proof in~\cite{Asmussen1998} relied on the local asymptotics for $\pr(M\in(x,x+T) $ found in~\cite{BD1996} and, independently, in~\cite{AKKKT}. 
Foss and Zachary~\cite{FZ2003} pointed out the necessity of the condition $F\in\mathcal S^*$ and  extended~\eqref{asym.Mtau} to the 
case of an arbitrary stopping time $\sigma$ with the finite mean 
$\mathbf E\sigma<\infty$. 
Then, Foss, Palmowsky and Zachary~\cite{FPZ2005}  found 
the asymptotics $\pr(M_\sigma>x)$, 
for a more general class of stopping times  $\sigma$,  including those  that  
may take infinite values or have infinite mean. 
They also proved that these
asymptotics hold uniformly in all stopping times. 
A short proof of~(\ref{asym.Mtau}) may be found in~\cite{D2004},\cite{D2005} and~\cite{DW2012}. 
The former proof relies on the local asymptotics for $\pr(M_\tau\in(x,x+T])$ and the latter proof uses the martingale properties of $\pr(M>x)$.
The local asymptotics for $\pr(M_\tau\in (x,x+T])$ were found in~\cite{DS2007}.

We will now introduce several subclasses of heavy-tailed distributions that will be used in the text. 
\begin{definition}
A distribution function $F$ is {\it (right) long tailed} ($F\in\mathcal L$) 
if, for any fixed $y>0$,  
\[
\mathbf P(\xi>x+y\mid\xi>x)=\frac{\overline F(x+y)}{\overline F(x)}\to 1 
, \quad x\to\infty.
\]
\end{definition}

An important subclass of heavy-tailed distributions is 
a class of subexponential distributions introduced independently 
by Chistyakov~\cite{Chist} and Chover {\it et al}~\cite{CNW}.
\begin{definition}
A distribution function $F$ on $\mathbf R^+$ is {\it subexponential} 
($F\in\mathcal S$) if 
 $\overline F(x)>0$ for all $x$ and 
\begin{equation}\label{sub1}
\frac{\mathbf P(\xi_1+\xi_2>x)}{\mathbf P(\xi_1>x)}
=\frac{\overline F^{2*}(x)}{\overline F(x)}\to 2, \quad \mbox{ as } x\to\infty,
\end{equation}
where $\xi_1,\xi_2$ are independent random variables with a common 
distribution function $F$.
\end{definition}

Sufficient conditions for a distribution to
belong to the class~$\mathcal S\,$ may be found,
or example, in~\cite{Chist},~\cite{Kl},~\cite{Pit} and~\cite{Teugels}.
The class~$\mathcal S\,$ includes, in particular,
the following distributions on $[0,\infty)$:
\begin{itemize}
    \item[(i)]~the Pareto distribution with the tail 
$\overline G(x)=(\frac{\kappa}{\kappa+x})^\alpha$, where $\alpha,\kappa>0$;
\item[(ii)]~the lognormal distribution with the density
$e^{-(\ln x-\ln\alpha)^2/2\sigma^2}/x\sqrt{2\pi\sigma^2}$
with ${\alpha>0}$;
\item[(iii)]~the Weibull distribution with the tail
$\overline G(x)=e^{-x^\alpha}$ with ${\alpha\in(0,1)}$.
\end{itemize}
Another subclass of heavy-tailed distribution 
is a class of distributions with dominated
varying tail. 
\begin{definition}
A distribution function $F$ is a {\it dominated varying} tail  distribution function 
 (${F\in\mathcal D}$)  if 
\[
\sup_{x>0}\frac{\overline F(x/2)}{\overline F(x)}<\infty.
\]
\end{definition}
A distribution function from $\mathcal D$ is not always subexponential. 
Indeed, all subexponential distributions are long-tailed, but there are some dominated
varying distributions which are not long-tailed, see~\cite{EG1}  and~\cite{Goldie} for a counterexample. 
 However, Kl\"uppelberg~\cite{Kl} proved that if the mean $\int_0^\infty
 \overline F(y) dy$ is finite then 
$\mathcal L\cap \mathcal D\subset \mathcal S^*\subset \mathcal S$. All regularly varying 
distribution functions belong to $\mathcal D$.
\begin{definition}
A distribution function $F$ 
is {\it regularly varying} with index $-\alpha$ if, for all $\lambda>0$, 
\[
\frac{\overline F(\lambda x)}{\overline F(x)}\to \lambda^{-\alpha},\quad
x\to\infty.
\]
\end{definition}
Examples of regularly varying distribution functions are the Pareto distribution
function and $G$ with the tail $\overline G(x)\sim 1/x^\alpha\ln^\beta x$
An extensive survey of the regularly varying distributions may be found in~\cite{BGT}.
It is shown in~\cite{Chist} that any subexponential distribution is long-tailed
with necessity. 
The converse is not true, see~\cite{EG1} for a counterexample.

When the mean is finite, the  derivation of 
(\ref{asym.Mtau}) in~\cite{Asmussen1998},\cite{FZ2003} and\cite{D2005}   
heavily relied on the local asymptotics 
$\mathbf P(M\in (x,x+c])$ for a fixed $c>0$ as $x\to \infty$. 
In the infinite-mean case,  these local asymptotics are not known. 
It seems that it can be found only in some particular cases. The reason for that  are complications in the local renewal  theorem in the infinite mean case, 
see~\cite{CaravennaDoney2016} for the complete solution of the local renewal  in the infinite mean case and its history.  
 Therefore, we propose a slightly different approach: it appears that it is sufficient to prove directly that 
\begin{equation}\label{weak}
\int_0^x\pi(du)\overline F(x-u)\sim \overline F(x).
\end{equation}

For that, we use a introduce a new class $\mathcal S_F$ of heavy-tailed distributions,
\begin{definition}\label{def1}
Let $F$ be a distribution function. 
A distribution function $G$ on ${\bf R}^+$ belongs to ${\mathcal
S}_F$   ($G\in {\mathcal S}_F$) if
\begin{equation}\label{defIG}
\int_0^x G(du)\overline F(x-u)\sim\overline F(x).
\end{equation}
\end{definition}
This class is a natural extension of the class of subexponential distributions. 
Indeed, it follows from the definition that 
$F$ is subexponential if and only if 
$F\in\mathcal S_F$. Then we study properties of this class. 
These properties (as well as its proofs) are
rather close to that of subexponential distributions. 
Let $G_1$ be a distribution function on $\mathbb R^+$ with the distribution tail 
\[
    \overline G_1(x) = \frac{1}{K}\int_x^{\infty} \frac{t-x}{m(t-x)} F(dt),\quad x\ge 0. 
\]
The following theorems are the main results of this paper. 
\begin{theorem}\label{th3}
Suppose ${\bf E}\xi^-=\infty$ and condition {\rm(\ref{finite.M})} holds.
 If the distribution function $G_1\in\mathcal S_F$
then asymptotics {\rm(\ref{asym.Mtau})} hold, i.~e.
\[
\mathbf P (M_\tau>x)\sim \mathbf E\tau\overline F(x).
\] 
\end{theorem} 

 \begin{theorem}\label{corr1}
Let $\mathbf E\xi_1^-=\infty$ and either of the following conditions holds 
\begin{enumerate}
\item[{\rm (a)}]  $F\in{\mathcal S}^*$; 
\item[{\rm (b)}] $F\in{\mathcal L}\cap{\mathcal D}$ and condition 
{\rm(\ref{finite.M})} holds. 
\end{enumerate}
Then $G_1\in\mathcal S_F$.
\end{theorem}

\section{Class ${\mathcal S}_F$ and its basic properties}\label{s2}

Definition~\ref{def1} may be rephrased as follows.
Consider independent random variables $\psi\ge 0$ and $\xi$ 
with distributions $G$ and $F$ respectively. 
Then $G\in {\mathcal S}_F$ if and only if 
\[ 
    \pr(\xi+\psi>x, \psi\le x)\sim \pr(\xi>x). 
\]

 Basic properties of the class ${\mathcal S}_F$ are very close to those of
the class of subexponential distributions(see Lemmas~\ref{prop}--\ref{lem3}). 
For a fine account of the theory of subexponential and local subexponential distributions we refer to~\cite{AFK} and ~\cite{FKZ2013}.

\begin{lemma}\label{prop} Let $G$ be a distribution function on ${\mathbf R}^+$. 
Then $G\in {\mathcal S}_F$ if and only
if there exists a function $h(x)\uparrow\infty, h(x)<x/2$ such that:
\begin{enumerate}
    \item[(1)]\label{c1} $\overline F(x-h(x))\sim \overline F(x)$;
    \item[(2)]\label{c2} $G(x-h(x),x]=o(\overline F(x))$;
    \item[(3)]\label{c3}
$
 \int_{h(x)}^{x-h(x)} G(du)\overline F(x-u)= o(\overline F(x)).
$
\end{enumerate}
\end{lemma}
{\sc Proof of Lemma \ref{prop}.} 
It is clear that if such a function $h(x)$ exists, then $F\in\mathcal L$ and
$G(x-t,x]=o(\overline F(x))$ for all fixed $t>0.$ Conversely, if
$F\in\mathcal L$ and $G(x-t,x]=o(\overline F(x))$ for some fixed $t>0$,
then one can construct a function $h(x)$ satisfying conditions (\ref{c1})
and (\ref{c2}).

First, assume $G\in S_F$. Fix any $t>0$. Then, 
\begin{eqnarray*}
\int_0^x G(du)\overline F(x-u)&=&
\left(\int_0^t+\int_t^{x-t}+\int_{x-t}^{x}\right) G(du)\overline F(x-u)\\
&\ge& G[0,t]\overline F(x)+G(t,x-t]\overline F(x-t)+G(x-t,x]\overline F(t).
\end{eqnarray*}
By dividing both sides by $\overline F(x)$, letting $x$ to infinity and rearranging the terms, we obtain 
\[
1\ge\lim\sup_{x\to\infty}\frac{\overline F(x-t)}{\overline
F(x)}\ge \lim\inf_{x\to\infty}\frac{\overline F(x-t)}{\overline
F(x)}\ge 1
 \]
 and 
 \[
 \limsup_{x\to\infty}\frac{G(x-t,x]}{\overline F(x)}=0
 \]
 Consequently, there exists a function
 $h(x)$ such that (\ref{c1}) and (\ref{c2}) hold. For this function, 
 condition (\ref{c3}) holds. 

Conversely, suppose that there is a function $h(x)$ satisfying 
(\ref{c1})--(\ref{c3}). Condition (\ref{c1}) implies 
\[
 \int_0^{h(x)}G(dy)\overline F(x-y)\sim \overline F(x)
 \int_0^{h(x)}G(dy)\sim \overline F(x),
 \]
and condition (\ref{c2}) implies 
\[
0\le\int_{x-h(x)}^{x}G(dy)\overline
F(x-y)\le 
G(x-h(x),x]=o(\overline F(x)).
\] 
Using condition
(\ref{c3}) we obtain the required result $G\in{\mathcal S}_F.$

 \qed

\begin{lemma}\label{conv}(convolution closure)
Let distribution functions $G_1, G_2$ belong to ${\mathcal S}_F$. 
Then  $G_1*G_2\in {\mathcal S}_F.$
\end{lemma}
{\sc Proof.}
Take a function $h(x)$ satisfying conditions (\ref{c1})--(\ref{c3}) 
of Lemma \ref{prop}  
 for both distributions $G_1$ and $G_2$ simultaneously. Then,
\begin{eqnarray*}
&\int_0^x G_1*G_2(du)\overline F(x-u)=
\int_0^{x} G_1(du)\int_0^{x-u}G_2(dv)\overline F(x-u-v)
&\\
&=
\left(\int_0^{x-h(x)}+\int_{x-h(x)}^{x}\right) G_1(du)\int_0^{x-u}G_2(dv)\overline F(x-u-v)
\equiv
I_1(x)+I_2(x).&
\end{eqnarray*}
By conditions~(\ref{c1}) and~(\ref{c3}) of Lemma~\ref{prop}, 
\begin{align*} 
    I_1(x)&=\int_0^{x-h(x)} G_1(du)\int_0^{x-u} G_2(dv)\overline F(x-u-v)\\
          &\sim \int_0^{x-h(x)} G_1(du)\overline F(x-u)\sim\overline F(x) 
\end{align*}
and, by condition (\ref{c2}), 
\[ 
    I_2(x)\le 
G_1(x-h(x),x]= o(\overline F(x)). 
\]
\qed

By induction, Lemma~\ref{conv} yields 
\begin{corollary}\label{subex}
 Let $G\in {\mathcal S}_F$. Then $G^{*n}\in {\mathcal S}_{F}$ for any $n\ge 1$.
\end{corollary}

Throughout, $G^{*0}$ denotes a distribution degenerated at $0$.

\begin{lemma} \label{majorant}
Let $G\in{\mathcal S}_F$. Then, for any $\varepsilon>0$, there exists
$A\equiv A(\varepsilon)>0$ such that, for any integer $n\ge 0$ and for any $x\ge 0$,
\begin{eqnarray*}
\int_0^x G^{*n}(du)\overline F(x-u) &\le& A(\varepsilon)(1+\varepsilon)^n
\overline F(x).
\end{eqnarray*}
\end{lemma}

{\sc Proof.}
Take any $\varepsilon>0$. Since $G\in{\mathcal S}_F$, there exists
$x_0>0$ such that 
 \begin{equation}\label{1e} 
\int_0^x G(du)\overline
F(x-u) \le(1+\varepsilon)\overline F(x), \quad x\ge x_0.
\end{equation} 
Put $A\equiv \frac{1}{\overline F(x_0)}$. We  use  
induction arguments.  For $n=0$ the assertion clearly holds.
 Suppose that the
assertion is true for $n-1$ and prove it for $n$. For $x<x_0$, 
\[
\frac{\int_0^x G^{*n}(du)\overline
F(x-u)}{\overline F(x)}\le \frac{1}{\overline F(x)}\le
\frac{1}{\overline F(x_0)}=A \le A(1+\varepsilon)^n 
\] 
Further, for $x\ge x_0$,
\begin{eqnarray*}
\displaystyle\int_0^x G^{*n}(du)\overline F(x-u) &=& \int_0^ x
G(dy) \int _0^{x-y}G^{*(n-1)}(dv)\overline F(x-y-v)\nonumber\\
&\displaystyle\le & A(1+\varepsilon)^{n-1}\int_0^{x}G(dy) \overline
F(x-y)\le A(1+\varepsilon)^n\overline F(x).
\end{eqnarray*}
The latter follows from (\ref{1e}).
\qed

Let $\{\zeta_n\}$ be a sequence of i.i.d. non-negative random variables with 
a common distribution $G$, and let   $\nu$ be a
 random stopping time independent of $\{\zeta_n\}$.
 Put $X_n=\zeta_1+\cdots+\zeta_n$. Then the
distribution of the randomly stopped sum $X_\nu$ is
\[
G_\nu(x)\equiv\mathbf P(X_\nu\le x)=\sum_{n\ge 0} \mathbf P(\nu=n)G^{*n}(x).
\]
\begin{lemma}\label{lem3}
Let $G$ belong to ${\mathcal S}_F$. 
Assume that  $\e (1+\delta)^{\nu}<\infty$ for some $\delta>0$. 
Then $G_\nu$ belongs to ${\mathcal S}_F$.
\end{lemma}
{\sc Proof.} The result follows from Corollary~\ref{subex}, Lemma~\ref{majorant} and from the dominated convergence theorem.

\qed

It is known~\cite[Theorem 3.2]{Kl} that 
if $F\in{\mathcal S}^*$ then $F$ is subexponential, i.e. $F\in\mathcal S_F$. 
In the following Lemma, we generalize this assertion 
to obtain sufficient conditions for $G\in {\mathcal
S }_F$. Another extension may be found in~\cite[Lemma 9]{DFK}.

\begin{lemma}\label{ssstar}
Let $F\in {\mathcal S}^*$ and $G$ on $\mathbf R^+$ be such that
$G(x-1,x]=o(\overline F(x))$. Then $G\in {\mathcal S}_F$.
\end{lemma}
{\sc Proof.} It follows from  $F\in {\mathcal S}^*$ that $F\in {\mathcal L}$, see~\cite[Theorem 3.2]{Kl}. Then there exists a function $h(x)$ satisfying
conditions (\ref{c1}) and (\ref{c2}) of Lemma~\ref{prop}.
 Further, 
 \[ 
     \int_{h(x)}^{x-h(x)} G(dy)\overline
F(x-y)\le \sum_{k=[h(x)]-1}^{[x-h(x)]} G(k-1,k]\overline F(x-k), 
\] 
Since $G(x-1,x]=o(\overline F(x))$, for $x$ large enough 
$G(x-1,x]\le\overline F(x)$.
  Then,
  \[
  \sum_{k=[h(x)]-1}^{[x-h(x)]} G(k-1,k]\overline F(x-k)\le
\sum_{k=[h(x)]-1}^{[x-h(x)]} \overline F(k)\overline F(x-k-1)=o(\overline F(x)),
 \]
the latter is due to $F\in\mathcal S^*$. Then, condition $(\ref{c3})$  
of Lemma~\ref{prop} is satisfied and $G~\in~\mathcal S_F.$
  \qed
  
As a simple corollary of Lemma~\ref{ssstar}, 
we can obtain sufficient conditions for the convolution $F*G$ 
to belong to the class $\mathcal S^*\subset\mathcal S$.  

\begin{corollary}\label{convol}
Let $F$ and $G$ on $\mathbf R^+$ be such that 
$F\in {\mathcal S}^*$ and $G(x-1,x]=o(\overline F(x))$. 
Then $F*G$ belongs to $\mathcal S^*$.
\end{corollary}

Let $H$ be a non-decreasing function on $\mathbf R^+$ such 
that integral 
\[
\int_0^\infty H(dt)\overline F(t) \mbox{ is finite.}
\]
We assume  that   $H$ is subadditive, i.e. if $H(x+y)\le H(x)+H(y)$ for all $x,y\ge 0$. 
Let $H(x,y]=H(y)-H(x)$. 
Consider a distribution function $G_H$ on $\mathbf R^+$ with the tail
distribution  
\begin{equation}\label{gu}
\overline G_H(x)\equiv\min
\Biggl(1,\int_0^\infty \overline F(t+x)\,H(dt)\Biggr),\quad x\ge 0.
\end{equation}
Integrating (\ref{gu}) by parts, we obtain an equivalent 
representation for $G_H$,
\begin{eqnarray}\label{gu2}
\overline G_H(x) &=& \min
\Biggl(1,\int_x^\infty H(0,t-x]\,F(dt)\Biggr),\quad x\ge 0.
\end{eqnarray}
We now establish some properties of $G_H$, which will be used in  the next Section.
\begin{lemma}\label{osmall}
Let $F\in\mathcal L$ and $H(x-1,x]\to 0$ as $x\to\infty$. Then,
\[
 G_H(x-1,x]=o(\overline F(x)).
\]
\end{lemma}
{\sc Proof.} 
It follows from definition that, for all sufficiently large $x$,
\[ 
    G_H(x-1,x]=
\int_{x-1}^\infty  H((t-x)^+,t-x+1] F(dt) 
\] 
Since $F\in\mathcal L$, there exists a function $h(x)$ such that 
$\overline F(x)\sim\overline F(x+h(x))$. Then, 
\[
\int_{x-1}^{x+h(x)} H((t-x)^+,t-x+1] F(dt)\le H(1) F(x-1,x+h(x)]
=o(\overline F(x))
 \]
and 
\[
\int_{x+h(x)}^\infty H(t-x,t-x+1] F(dt)\le
\sup_{y\ge h(x)}H(y,y+1] \overline F(x+h(x))=o(\overline F(x)).
\]
\qed

\begin{lemma}\label{cor3}
Let $F\in\mathcal L$, and let $H_1,H_2$ be subadditive functions such that 
$H_1(x-1,x]\to 0, H_2(x-1,x]\to 0$ and 
\begin{equation}\label{33}
0<\liminf \frac{H_1(x)}{H_2(x)}\le\limsup \frac{H_1(x)}{H_2(x)}<\infty.
\end{equation}
Then $$G_{H_1}\in\mathcal S_F\iff G_{H_2}\in\mathcal S_F.$$ 
\end{lemma}
{\sc Proof.} 
By Lemma~\ref{osmall}, $G_{H_1}(x-1,x]=o(\overline F(x))$
and  $G_{H_2}(x-1,x]=o(\overline F(x))$. Therefore, there exists a
function $h(x)$ such that conditions (\ref{c1}) and (\ref{c2}) of
Lemma~\ref{prop} hold for both distribution functions $G_{H_1}$ and $G_{H_2}$.
 Integrating by parts (\ref{c3}) and using 
(\ref{c1}) and (\ref{c2}) , we obtain 
\[
\int_{h(x)}^{x-h(x)}G_{H_i}(dy)\overline F(x-y)=
\int_{h(x)}^{x-h(x)}F(dy)\overline G_{H_i}(x-y,x]+o(\overline F(x)), i=1,2.
\] 
For any $y<x$, we have 
\begin{eqnarray*} 
\overline G_{H_i}(x-y,x]=
\int_{x-y}^x  H_i(0,t-x+y] F(dt) 
+\int_{x}^\infty  H_i(t-x,t-x+y] F(dt).
\end{eqnarray*}
Due to the subadditive property, 
\[
\int_{h(x)}^{x-h(x)}F(dy)\int_{x}^\infty  H_i(t-x,t-x+y]
F(dt)\le\int_{h(x)}^{x-h(x)}F(dy) H_i(y)\overline F(x)=o(\overline F(x)). 
\]  
Hence condition (\ref{c3}) of Lemma \ref{prop} 
holds for distribution functions $G_{H_1}$ and $G_{H_2}$ if and only if  
\begin{equation}\label{34}
\int_{h(x)}^{x-h(x)}F(dy)\int_{x-y}^x  H_i(0,t-x+y] F(dt)=o(\overline F(x)).
\end{equation}
Then the assertion of Lemma follows from (\ref{33}).
\qed

\begin{lemma}\label{lem12}
Assume that  $F\in\mathcal L\cap \mathcal D$ and 
$H(x-1,x]\to 0$ as $x\to\infty$. Then $G_H\in \mathcal S_F.$ 
\end{lemma}
{\sc Proof.}
It follows from Lemma \ref{osmall} that
$G_H(x-1,x]=o(\overline F(x))$. Therefore, there exists a function satisfying
conditions (\ref{c1}) and (\ref{c2}). From the proof of Lemma
\ref{cor3}, it is clear that if (\ref{34}) holds then $G_H\in\mathcal S_F$. 
Using the subadditive property of $H$, we
obtain 
\[
\int_{h(x)}^{x/2}F(dy)\int_{x-y}^x  H(0,t-x+y] F(dt)\le
 \int_{h(x)}^{x/2}F(dy)H(0,y]  F(x/2,x]=o(1) \overline F(x/2).
\]
 Further, 
 \[
\int_{x/2}^{x-h(x)}F(dy)\int_{x-y}^x  H(0,t-x+y] F(dt)\le
\overline F(x/2) \int_{h(x)}^\infty F(dt)H(0,t]=o(1) \overline F(x/2).
\]
It follows from $F\in\mathcal D$ that $o(1) \overline
F(x/2)=o(\overline F(x))$. 
\qed


\section{Proof of the main results}

First recall a well-known construction of ladder moments and ladder heights~\cite[Chapter XII]{F}. 
Let 
\[
\eta = \min\{n\ge 1: S_n>0\}\le \infty
\] 
be the first (strict) ascending ladder epoch and put 
\[
p=\mathbf P\{\eta=\infty\}=\mathbf P (M=0). 
\] 
Let
$\{\psi_n\}_{n\ge 1}$ be a sequence of i.i.d.r.v.'s distributed as
\[ 
\mathbf P(\psi_1\in B)\equiv G_+(B)=\mathbf P(S_\eta\in B |
\eta<\infty). 
\]
Let $\nu$ be a random variable, independent of
the above sequence, such that 
$\mathbf P (\nu=n)=p(1-p)^n, n=0,1,2,\ldots $ 
Then
\begin{equation}\label{Msum}
M\stackrel{d}=\psi_1+\ldots+\psi_\nu.
\end{equation}
We start with proving an auxiliary assertion. 
\begin{lemma}\label{th1}
Let $G_+\in {\mathcal S}_F$. Then asymptotics {\rm(\ref{asym.Mtau})} hold, i.~e.
$$ \mathbf P (M_\tau>x)\sim \mathbf
E\tau\overline F(x). $$
\end{lemma}
{\sc Proof.} The proof  is carried out via standard arguments: we obtain the lower and the upper bounds, which are asymptotically equivalent. 
Let us start with the lower bound, which is valid without any
assumptions on $F$ and $G_+$. Fix a positive integer
$N$, then for any $x>0$,
\begin{eqnarray*}
 \mathbf P(M_\tau>x)&\ge&\sum_{n=0}^N 
 \mathbf P(\tau>n, \max_{0\le i\le n} S_i\le x, S_{n+1}>x)\\
&\ge& 
\sum_{n=0}^N 
 \mathbf P(\tau>n, \max_{0\le i\le n} S_i\le x, \xi_{n+1}>x)\\
 &=&(1+o(1))\overline F(x)
 \sum_{n=0}^N \mathbf P(\tau>n).
\end{eqnarray*}
Let $N\to \infty$ to obtain 
\[
\mathbf P(M_\tau>x)\ge (1+o(1))\mathbf E\tau\overline F(x).
\]

Now turn to the upper bound. 
For any $x\geq 0$,
\begin{eqnarray*}
\displaystyle\mathbf P(M_\tau>x)&\le & 
\sum_{n=0}^{\infty}
 \mathbf P(\tau>n,  S_n\le x, S_{n+1}>x)\\&=&\int_0^x\sum_{n=0}^{\infty}
 \mathbf P(S_1>0,\ldots, S_n>0,  S_n \in du)\overline F(x-u)\\
&\displaystyle=&\mathbf E\tau\int_0^x \mathbf P(M\in du)\overline F(x-u),
\end {eqnarray*}
for the latter see \cite[Chapter XII, (2.7)]{F}. 
Finally, it follows from $G_+\in{\mathcal S}_F$, relation (\ref{Msum}) and 
Lemma~\ref{lem3} that the distribution function of $M$ 
belongs to $\mathcal S_F$, that is
$$ \int_0^x \mathbf P(M\in du)\overline F(x-u)\sim
\overline F(x). $$ \qed


Now let us introduce a few more definitions. Let $\chi =-S_{\tau}$ be the
absolute value of the first non-positive sum and let  $ G_-(x)\equiv\mathbf P(\chi\le x) $ be its distribution function. 
Define a renewal function 
$$
H_-(x) =  \sum_{n=0}^\infty G_-^{*n}(x), \quad x\ge 0. 
$$
Then $\psi_1$ is distributed as follows \cite[Chapter XII]{F}:
\begin{eqnarray}\label{psiH}
\mathbf P(\psi_1>x)
\equiv  \overline G_+(x)
&=&\frac{1}{1-p}
\int_0^\infty \overline F(u+x) H_-(du)
\end{eqnarray}

We will need the following asymptotic estimates for  the renewal function $H$ 
\begin{proposition}(see~\cite[Corollary 2]{DFK}) 
    Suppose $\e[\xi_1^-]=\infty$ and condition~\eqref{finite.M} holds. Then, 
\begin{eqnarray}\label{p1}
p\le\liminf_{x\to\infty} \frac{H_-(x)m(x)}{x}\le\limsup_{x\to\infty} 
\frac{H_-(x)m(x)}{x}\le 2p
\end{eqnarray}
\end{proposition}
{\sc Proof of Theorem~\ref{th3}.} We will  prove that 
$G_1\in\mathcal S_F$ implies that $G_+\in\mathcal S_F$. 
For that, we verify the conditions of Lemma~\ref{cor3}. 

First, by the properties of renewal functions, 
$H_-$ is subadditive. It
follows from $\mathbf E\xi_1^-=\infty$ that $\mathbf E\chi=\infty.$ 
This,  the Key Renewal Theorem  implies that $H(x-1,x]\equiv H(x)-H(x-1)\to 0, x\to\infty$.

Second, the function $x/m(x)$ is
 subadditive as well: 
 \[
\frac{x+y}{m(x+y)}=\frac{x}{m(x+y)}+\frac{y}{m(x+y)}\le 
\frac{x}{m(x)}+\frac{y}{m(y)}.
\]
The latter inequality holds  since $m(x)$ is non-decreasing. 
Further,  the function $x/m(x)$ 
is non-decreasing, since 
\[
\frac{d}{dx}\frac{x}{m(x)}=\frac{m(x)-xm'(x)}{m^2(x)}=\frac{m(x)-x\pr(\xi_1^->x)}{m^2(x)}\ge 0. 
\]
 Hence, 
 \[
0\le\frac{x}{m(x)}-\frac{x-1}{m(x-1)}\le\frac{1}{m(x-1)}\to 0.
\]
Then,   Lemma~\ref{cor3} and (\ref{p1}) imply that 
 $G_+\in\mathcal S_F$ if and
 only if $G_1\in\mathcal S_F.$ 
Consequently, $G_+\in\mathcal S_F$ and, by Lemma~\ref{th1}, 
asymptotics (\ref{asym.Mtau}) hold.  
\qed 

{\sc Proof of Theorem~\ref{corr1}.}
Sufficiency of (a) follows from  Lemma~\ref{ssstar} and Lemma~\ref{osmall}. 
Sufficiency of (b) follows from Lemma~\ref{lem12}. 
\qed

{\bf Acknowledgements} 

The author is most grateful to Sergey Foss for drawing attention to this problem, valuable comments regarding the manuscript  and a number of useful discussion. 


\end{document}